\documentclass[11pt]{article}
\usepackage{amssymb}

\def\OV#1{\overline {#1}}
\def\hang{\hangindent\parindent}
\def\textindent#1{\indent\llap{#1\enspace}\ignorespaces}
\def\re{\par\hang\textindent}

\def\QED{\hfill{$\Box$}}
\def \r{\rightarrow}

\def\mapright#1#2{\smash{\mathop{\longrightarrow}\limits^{#1}_{#2}}}

\def\v5{\vskip .5truecm}

\def\LM{{\bf LM}}\def\LT{{\bf LT}}
\def\B{{\cal B}} \def\LC{{\bf LC}} \def\G{{\cal G}} \def\FRAC#1#2{\displaystyle{\frac{#1}{#2}}}
\def\SUM^#1_#2{\displaystyle{\sum^{#1}_{#2}}} \def\J{{\bf J}} \def\LR{\langle X\rangle}
\def\LH{{\bf LH}}
\def\FRAC#1#2{\displaystyle{\frac{#1}{#2}}}\def\GFA{G(A)} \def\GFM{G(M)}
\def\SUM^#1_#2{\displaystyle{\sum^{#1}_{#2}}}\def\G{{\cal G}}

\title{Noncommutative Gr\"obner Bases for \\ Almost Commutative Algebras\thanks{Project supported by
the National Natural Science Foundation of China (10571038).}}
\author{Huishi Li\\
{\small Department of Applied Mathematics}\\
{\small College of Information Science and Technology}\\
{\small Hainan University}\\
{\small  Haikou 570228, China}}
\date{}

\begin{document}
\maketitle
\begin{center}
\begin{minipage}{120mm}
{\small {\bf Abstract.} Let $K$ be an infinite field and $K\langle
X\rangle =K\langle X_1,...,X_n\rangle$ the free associative algebra
generated by $X=\{ X_1,...,X_n\}$ over $K$. It is proved that if $I$
is a two-sided ideal of $K\langle X\rangle$ such that the
$K$-algebra $A=K\langle X\rangle /I$ is almost commutative in the
sense of [3], namely, with respect to its standard
$\mathbb{N}$-filtration $FA$, the associated $\mathbb{N}$-graded
algebra $G(A)$ is commutative, then $I$ is generated by a finite
Gr\"obner basis. Therefor, every quotient algebra of the enveloping
algebra $U(\mathbf{g})$ of a finite dimensional $K$-Lie algebra
$\mathbf{g}$ is, as a noncommutative algebra of the form $A=K\langle
X\rangle /I$, defined by a finite Gr\"obner basis in $K\langle
X\rangle$.}
\end{minipage}\end{center}{\parindent=0pt\vskip 6pt
{\bf 2000 Mathematics Classification} Primary 16W70; Secondary 16Z05.\vskip 6pt
{\bf Key words} Almost commutative algebra, Filtration, gradation, Gr\"obner basis}

\vskip 1truecm
\def\GR{Gr\"obner} \def\B{{\cal B}}
\def\GL{\ge_{grlex}} \def\LM{{\bf LM}} \def\LC{{\bf LC}} \def\O{{\bf
O}} \def\J{{\bf J}}\def\U{{\bf U}}\def\LH{{\bf LH}}\def\LT{{\bf LT}}
\def\LR{\langle X\rangle}
\def\FRAC#1#2{\displaystyle{\frac{#1}{#2}}}\def\GFA{G(A)} \def\GFM{G(M)}
\def\SUM^#1_#2{\displaystyle{\sum^{#1}_{#2}}}\def\G{{\cal G}}
\section*{1. Introduction}
Let $K$ be an infinite field and let $K\LR =K\langle X_1,X_2,...,X_n\rangle$ be
the free associative $K$-algebra generated by $X=\{ X_1,X_2,...,X_n\}$ over $K$. It is
a well-known fact that even if a two-sided ideal $J$ of $K\LR$ is finitely generated,
$J$ does not necessarily have a finite Gr\"obner basis in the sense of [9]. However, it
was proved in [4] that if the quotient algebra $A=K\LR /J$ is commutative,
then, after a general linear change of variables (if necessary),
$J$ has a finite Gr\"obner basis in $K\LR$ (the construction of a Gr\"obner basis for $J$
is mentioned in the beginning of section 3 below). This, indeed,
makes another algorithmic way to study a commutative algebra via its noncommutative
Gr\"obner presentation, and its effectiveness  may be illustrated, for example, by the work
of [1], [2], [5] and [11]. In this note, we use the result of [4] and the filtered-graded
transfer of Gr\"obner bases [7] to show that for a two-sided ideal $I$ of
$K\langle X\rangle$, if the $K$-algebra $A=K\langle X\rangle /I$ is almost commutative
in the sense of [3], namely, with respect to its standard $\mathbb{N}$-filtration $FA$,
the associated $\mathbb{N}$-graded algebra $G(A)$ is commutative, then $I$ is generated
by a finite Gr\"obner basis. After reaching the main result in section 2, we discuss
in section 3 how this result may be realized computationally.\par
Throughout the note, we fix the infinite field $K$, the finite set $X=\{ X_1,X_2,...,X_n\}$
and the corresponding free $K$-algebra $K\LR =K\langle X_1,X_2,...,X_n\rangle$.
\v5

\section*{2. The Main Result}
Since the free $K$-algebra $K\LR =K\langle X_1,X_2,...,X_n\rangle$ has its standard
$K$-basis $\B$ consisting of all monomials (words)
$X_{j_1}X_{j_2}\cdots X_{j_s}$, $s\in\mathbb{N}$, to be convenient, we denote monomials
in $\B$ by $u,~v,~w$, $\cdots$.
Consider the $\mathbb{N}$-gradation $K\LR =\oplus_{p\in\mathbb {N}}K\LR_p$ of $K\LR$ with
$K\LR_p=K$-Span$\{ u\in\B~|~d(u)=p\}$, where $d(u)$ stands for the degree (length) of $u$.
If $I$ is a two-sided ideal of $K\LR$, then the $K$-algebra $A=K\LR /I$ has the {\it
standard $\mathbb{N}$}-{\it filtration}
$$F_0A=K\subset F_1A\subset\cdots\subset F_pA\subset\cdots$$
where for each $p\in\mathbb{N}$,
$$F_pA=(\oplus_{k\le p}K\LR_k+I)/I,$$
which defines the associated $\mathbb{N}$-graded $K$-algebra
$G(A)=\oplus_{p\in\mathbb{N}}G(A)_p$ of $A$ with $G(A)_p=F_pA/F_{p-1}A$. Recall from the
literature [3] that $A$ is called an {\it almost commutative algebra} if the
associated graded $K$-algebra $G(A)$ of $A$ is commutative. In [3] it was proved that
{\it $A$ is almost commutative if and only if $A$ is a homomorphic image of the universal
enveloping algebra $U(\mathbf{g})$ of some finite dimensional $K$-Lie algebra
$\mathbf{g}$.} For instance, the $n$-th Weyl algebra $A_n(K)$ is an almost commutative
algebra and it is a homomorphic image of the $2n+1$-dimensional Heisenberg Lie algebra.
So the class of almost commutative $K$-algebras consists of quotient algebras
of enveloping algebras of finite dimensional $K$-Lie algebras. It is well-known that the
study of quotient algebras of enveloping algebras (for example, the primeness, primitivity,
simplicity, etc) has been very important in the finite dimensional Lie-theory. The reader
is referred to ([10] Chapter 8 and Chapter 14) for details on this topic.
The result and its proof given below provide a way to study the structural properties of almost commutative $K$-algebras algorithmically, or
more precisely, to study such algebras by using both commutative and noncommutative Gr\"obner bases.
{\parindent=0pt\v5
{\bf 2.1. Theorem} If the $K$-algebra $A=K\LR /I$ is almost commutative, then  $I$ has a finite Gr\"obner basis. Equivalently,
for a finite dimensional $K$-Lie algebra $\mathbf{g}$, every quotient algebra of the enveloping algebra
$U(\mathbf{g})$ of $\mathbf{g}$, viewed as a quotient of some free $K$-algebra,
is defined by a finite Gr\"obner basis.
\vskip 6pt
{\bf Proof} For an element $f\in K\LR$, if $f=F_p+F_{p-1}+\cdots +F_0$ with $F_i\in K\LR_i$
and $F_p\ne 0$, then we write $d(f)=p$ for the degree of $f$ and $\LH (f)=F_p$ for the leading homogeneous part of $f$. Consider
the standard $\mathbb{N}$-filtration $FA$ of $A$ and its associated $\mathbb{N}$-graded
algebra $G(A)$, as mentioned before. By ([7] Chapter III Proposition 3.1, or [8]),
$G(A)\cong K\LR /\langle\LH (I)\rangle$, where $\langle\LH (I)\rangle$ is the graded two-sided ideal of
$K\LR$ generated by $\LH (I)=\{ \LH (f)~|~f\in I\}$. Since $A$ is almost commutative, $K\LR /\langle\LH (I)\rangle$ is
a commutative $K$-algebra.}\par
Now, let $K[\mathbf{x}]=K[x_1,x_2,...,x_n]$ be the commutative polynomial $K$-algebra in variables $x_1,x_2,...,x_n$ and
consider the natural $K$-algebra homomorphism $\gamma$: $K\LR \r
K[\mathbf{x}]$ with $\gamma (X_i)=x_i$. Then $G(A)\cong K\LR /\langle\LH (I)\rangle \cong
K[\mathbf{x}]/\gamma (\LH (I))$ and by ([4] Theorem 2.1, Corollary 1.1), a finite Gr\"obner basis of
$\langle\LH (I)\rangle$ may be obtained by using a finite Gr\"obner basis of $\gamma (\langle\LH (I)\rangle )$.
Our aim is to lift the obtained Gr\"obner basis of $\langle\LH (I)\rangle$ to a
finite Gr\"obner basis of $I$ as described in ([7] Chapter III section 3). To this end, we need to use a graded monomial ordering.
Note that  $\gamma$ does not change
the degree (length) of monomials, that is, for any two monomials $u$, $v$ in the
standard $K$-basis $\B$ of $K\LR$,
$$d(u)<d(v)~\hbox{if and only if}~d(\gamma (u))<d(\gamma (v)).$$
It turns out that if we fix a certain graded monomial ordering $\prec$ (for example, the graded
lexicographic ordering) on $K[\mathbf{x}]$, then the {\it lexicographic extension} $\prec_{et}$
of $\prec$ to $K\LR$ (in the sense of [4]), which is defined for $u,v\in\B$ by
$$u\prec_{et}v~\hbox{if}~\left\{\begin{array}{l} \gamma (u)\prec\gamma (v)~\hbox{or}\\
\gamma (u)=\gamma (v)~\hbox{and}~u~\hbox{is lexicographically smaller than}~v,\end{array}\right.$$
is a graded monomial ordering. Thus, it follows from ([4] Theorem 2.1, Corollary 1.1) that $\langle\LH (I)\rangle$
has a finite Gr\"obner basis $G$ with respect to $\prec_{et}$. Since $\gamma (
\langle\LH (I)\rangle )$ is a graded ideal of $K[\mathbf{x}]$, by the construction of $G$
given in ([4] Theorem 2.1, or see the description given in the beginning of section 3 below),
we may assume that  $G =\{ G_1,G_2,...,G_s\}$ consists of
homogeneous elements in $\langle\LH (I)\rangle$. Suppose for $i=1,2,...,s$,
$$G_i=\sum H_i\LH (f_i)T_i,~\hbox{where}~H_i,~T_i~\hbox{are homogeneous elements and} ~f_i\in I.$$
If we put $f_i=\LH (f_i)+f_i'$, where $d(f_i')<d(f_i)$, then $\sum H_if_iT_i=g_i\in I$ and
$$g_i=\sum H_if_iT_i=\sum H_i\LH (f_i)T_i+\sum H_if_i'T_i=G_i+\sum H_if_i'T_i$$
has its leading homogeneous part $\LH (g_i)=G_i$, $i=1,2,...,s$. As $\prec_{et}$ is a graded
monomial ordering on $K\LR$, it follows from ([7] Chapter III Theorem 3.7) that
$\G =\{ g_1,g_2,...,g_s\}$ is a finite Gr\"obner basis of $I$ with respect to $\prec_{et}$. This completes the proof.\QED
\v5

\section*{3. Further Discussion}
Concerning the computational realization of Theorem 2.1, we first recall from [4]
how a noncommutative Gr\"obner basis may be constructed by using a commutative version.
Let $K\LR =K\langle X_1,X_2,...,X_n\rangle$, $K[\mathbf{x}]=K[x_1,x_2,...,x_n]$, and
$\gamma$: $K\LR\r K[\mathbf{x}]$ with $\gamma (X_i)=x_i$ be as before. The {\it lexicographic splitting} of $\gamma$
is defined as the $K$-linear map
$$\begin{array}{ccccc} \delta :&K[\mathbf{x}]&\mapright{}{}&K\LR&\\
&x_{i_1}x_{i_2}\cdots x_{i_r}&\mapsto&X_{i_1}X_{i_2}\cdots X_{i_r}&i_1\le i_2\le\cdots\le i_r\end{array}$$
For a monomial ideal $L$ in $K[\mathbf{x}]$, if $m=x_{i_1}\cdots x_{i_r}\in L$ with
$i_1\le\cdots\le i_r$, let ${\cal U}_L(m)$ denote the set of all monomials $u\in K[x_{i_1+1},...,x_{i_r-1}]$
such that neither $u\frac{m}{x_{i_1}}$ nor $u\frac{m}{x_{i_r}}$ lies in $L$. Following
([4] Theorem 2.1), if $\Gamma =\{ g_1,g_2,...,g_s\}$ is a minimal $\prec$-Gr\"obner basis for
an ideal $I\subset K[\mathbf{x}]$, then a minimal $\prec_{et}$-Gr\"obner basis for $J=\gamma^{-1}(I)\subset K\LR$
consists of $\{ X_jX_i-X_iX_j~|~1\le i<j\le n\}$ together with the elements $\delta (u\cdot g_i)$ for
each $g_i\in \Gamma$ and each $u\in {\cal U}_{\langle\LH (I)\rangle}(\LM (g_i))$, where
$\LM (g_i)$ stands for the $\prec$-leading monomial of $g_i$ and $\langle\LM (I)\rangle$ is the
monomial ideal generated by $\LM (I)=\{\LM (f)~|~f\in I\}$.
\v5
In view of the characterization of an almost commutative algebra [3], to realize Theorem 2.1 computationally,
we have to consider two coherent cases separately. In what follows, notations used in section 2 are
maintained.
{\parindent=0pt\v5
{\bf I.} $A=K\LR /I$, $I=\langle S\rangle$ and $G(A)$ is commutative.\par
Since $G(A)$ is commutative, it follows from the proof of Theorem 2.1 that we have
$$G(A)\cong K\LR /\langle\LH (I)\rangle\cong K[\mathbf{x}]/\gamma (\langle\LH (I)\rangle ).$$
Again by the proof of Theorem 2.1, in order to obtain a finite Gr\"obner basis of $I$
algorithmically, we need to have a finite Gr\"obner basis of $\gamma (\langle\LH (I)\rangle )$
so that we may use it to construct a Gr\"obner basis for $\langle\LH (I)\rangle$ as remarked above.
While this implies that we need first to know a generating set of $\gamma (\langle\LH (I)\rangle )$.
But the fact is that even if the generating set of $I$ is finite, say $S=\{ f_1,f_2,...,f_s\}$,
the equality $\langle\LH (I)\rangle =\langle\LH (f_1),\LH (f_2),...,\LH (f_s)\rangle$ is not
necessarily true (see [7] Chapter III). In other words, if we cannot find a generating set of
$\gamma (\langle\LH (I)\rangle )$ effectively, Theorem 2.1 is only theoretical. }\par
As an example, let us point out that if we know some (finite or infinite) Gr\"obner basis $\G =\{
g_i~|~i\in\Omega\}$ of $I$ and $\G$ contains all commutators $X_iX_j-X_jX_i$, $1\le i<j\le n$,
then $\langle\LH (I)\rangle =\langle\LH (g_i)~|~i\in\Omega\rangle$ by ([7] Chapter III).
{\parindent=0pt\v5
{\bf II.} $A=U(\mathbf{g})/I$ and $I=\langle S\rangle$.\par
Suppose $U(\mathbf{g})=K\LR/J$ with $J=\langle [X_i,X_j]-X_iX_j-X_jX_i~|~1\le i<j\le n\rangle$,
where $[X_i,X_j]$ is given by the Lie-brackets of generators of the finite dimensional $K$-Lie algebra
$\mathbf{g}$. Let $\OV I$ be the two-sided ideal of $K\LR$ such that $I=\OV I/J$. Then
$A=U(\mathbf{g})/I\cong K\LR /\OV I$. By the foregoing discussion, if we can find a
finite Gr\"obner basis of $\gamma (\langle\LH (\OV I)\rangle )$, then a finite Gr\"obner basis
for $\OV I$ may be constructed. }\par
Note that the standard filtration $FU(\mathbf{g})$ induces the standard filtration
$FA$. Consider the filtration $FI$ on $I$ induced by $FU(\mathbf{g})$ and its associated
graded ideal $G(I)$ in $G(A)$. Then it is well known that
$$G(A)=G(U(\mathbf{g})/I)\cong G(U(\mathbf{g}))/G(I)\cong K[\mathbf{x}]/\overline{G(I)},$$
where $\OV{G(I)}$ stands for the counterpart of $G(I)$ in $K[\mathbf{x}]$.
Thus, since both graded algebra epimorphisms
$$\begin{array}{l} \varphi :~k[\mathbf{x}]~\mapright{}{}~K[\mathbf{x}]/\OV{G(I)}\cong G(A)\\
\\
\psi :~k[\mathbf{x}]~\mapright{}{}~K[\mathbf{x}]/\gamma (\langle\LH (\OV I)\rangle )\cong G(A)\end{array}$$
agree on the generators $x_1,x_2,...,x_n$, we have $\OV{G(I)}=\gamma (\langle\LH (\OV I)\rangle )$. This
makes the chance for us to have a finite Gr\"obner basis of $\gamma (\langle\LH (\OV I)\rangle )$
by using the generating set $S$ of $I$. To see this, first recall that $U(\mathbf{g})$ is a solvable
polynomial algebra in the sense of [6]. Hence, starting with $S$, a noncommutative version of Buchberger Algorithm
produces a finite Gr\"obner basis $\G =\{ g_1,g_2,...,g_m\}$ for $I$. Furthermore, if the monomial ordering
used in producing $\G$ is a graded monomial ordering, then it follows from ([7] Chapter IV Theorem 2.1) that
$\sigma (\G )=\{\sigma(g_1),\sigma (g_2),...,\sigma (g_m)\}$ is a Gr\"obner basis for $G(I)$ with respect to the
same type of graded monomial ordering, where for each $g_i$, if $g_i\in F_pI-F_{p-1}I$, then
$\sigma (g_i)$ is the homogeneous element in $G(I)_p$ represented by $g_i$. Finally,
by passing to $\OV{G(I)}$ in $K[\mathbf{x}]$, we obtain a finite
Gr\"obner basis $\{ \OV g_1,\OV g_2,...,\OV g_m\}$ of $\gamma (\langle\LH (\OV I)\rangle )$. Note that
$U(\mathbf{g})=K\LR /J$, $I=\OV I/J$. Hence the preimage of $\G$ is contained in $\OV I$. Thus,
by the definition of $\delta (u\cdot\OV g_i)$ described in the beginning of this section,
the last step of the proof of Theorem 2.1 can be realized to give a Gr\"obner basis for $\OV I$.\par
In this case, the good news is that nowadays there has been the well-developed computer algebra system
SINGULAR: PLURAL (http://www.singular.uni-kl.de/index.html) which provides a programme named  \textsf{twostd} for computing a two-sided
Gr\"obner basis of a two-sided ideal in a solvable polynomial algebra.
\v5
We end this note by employing
an easy example to illustrate the procedure demonstrated in part I and II above.
Let $U(s\ell_2)$ be the enveloping algebra of the 3-dimensional
$K$-Lie algebra $s\ell_2=Ke\oplus Kf\oplus Kh$ subject to the relations
$[e ,f]=h,\quad [h ,e]=2e ,\quad [h ,f]=-2f$,
that is, $U(s\ell_2)=K[e,f,h]\cong K\langle X,Y,Z\rangle /P$ with $X\mapsto e$,
$Y\mapsto f$ and $Z\mapsto h$, where
$$P=\langle YX-XY+Z, ~ZX-XZ-2X, ~ZY-YZ+2Y\rangle .$$
Consider the two-sided ideal $I=\langle e^3,f^3,h^3-4h\rangle$ of $U(\mathbf{g})$, then
$U(\mathbf{g})/I\cong K\langle X,Y,Z\rangle /\OV I$, where
$$\OV I=\langle YX-XY+Z, ~ZX-XZ-2X, ~ZY-YZ+2Y,~X^3,~Y^3,~Z^3-4Z\rangle .$$
Let $\prec$ be the graded reverse lexicographic ordering on
$U(\mathbf{g})$ defined by  $e\prec f\prec h$. Then \textsf{twostd}
produces a two-sided Gr\"obner basis
$$\G =\left \{\begin{array}{l} e^3,~f^3,~h^3-4h,~eh^2+2eh,\\
fh^2-2fh,~2efh-h^2-2h,\\
~e^2f-eh-2e,~ef^2-fh,~e^2h+2e^2,~f^2h-2f^2\end{array}\right\}$$
for $I$ (Singular Manual A.6.1). Hence
$$\sigma (\G )=\left\{\begin{array}{l}
\sigma (e)^3,~\sigma (f)^3,~\sigma (h)^3,~\sigma (e)\sigma (h)^2,
~\sigma (f)\sigma (h)^2,\\
2\sigma (e)\sigma (f)\sigma (h),~\sigma (e)^2\sigma (f),~\sigma
(e)\sigma (f)^2, ~\sigma (e)^2\sigma (h),~\sigma (f)^2\sigma
(h)\end{array}\right\}$$ is a Gr\"obner basis of $G(I)$. If we use
the graded reverse lexicographic ordering  $x\prec y\prec z$ on the
polynomial $K$-algebra $K[x,y,z]$, then previous part II yields a
Gr\"obner basis
$$\{ x^3,~y^3,~z^3,~xz^2,~yz^2,~2xyz,~x^2y,~xy^2,~x^2z,~y^2z\}$$
for $\gamma (\langle\LH (\OV I)\rangle )$. By the construction of
$\delta (u\cdot g)$  described in the beginning of this section, it
may be checked directly that $\langle\LH (\OV I)\rangle$ has a
Gr\"obner basis
$$\left\{\begin{array}{l} YX-XY,~ZY-YZ,~ZX-XZ,~X^3,~Y^3,\\
Z^3,~XZ^2,~YZ^2,~2XYZ,~X^2Y,~XY^2,~X^2Z,~Y^2Z\end{array}\right\}.$$
Note that $\OV I$ contains the preimage of $\G$. It follows from the
last  step of the proof of Theorem 2.1 that $\OV I$ has the
Gr\"obner basis
$$\left \{\begin{array}{l} YX-XY,~ZY-YZ,~ZX-XZ,\\
X^3,~Y^3,~Z^3-4Z,~XZ^2+2XZ,\\
YZ^2-2YZ,~2XYZ-Z^2-2Z,\\
~X^2Y-XZ-2X,~XY^2-YZ,~X^2Z+2X^2,~Y^2Z-2Y^2\end{array}\right\} $$
with respect to $X\prec_{et}Y\prec_{et}Z$. \v5
\centerline{References}{\parindent=.65truecm\par \re{[1]} D. J.
Anick, On the homology of associative algebras, {\it Trans. Amer.
Math. Soc}., Vol. 296, 2(1986), 641--659. \re{[2]} D. Anick and
G.-C. Rota, Higher-order syzygies for the bracket algebra and for
the ring of coordinates of the Grassmannian, {\it Proc. Nat. Acad.
Sci. U.S.A.}, 88(1991), 8087¨C8090. \re{[3]} M. Duflo, Certaines
alg\`ebres de type finisont des alg\`ebres de Jacobson, {\it J.
Algebra}, 27(1973), 358--365. \re{[4]} D. Eisenbud, I. Peeva and B.
Sturmfels, Non-commutative Gr\"obner bases for commutative algebras,
{\it Proc. Amer. Math. Soc.}, 126(1998), 687-691. \re{[5]} E. Green
and R.Q. Huang, Projective resolutions of straightening closed
algebras generated by minors, {\it Adv. in Math.}, 110(1995),
314¨C333. \re{[6]} A.~Kandri-Rody and V.~Weispfenning,
Non-commutative Gr\"obner bases in algebras of solvable type, {\it
J. Symbolic Comput.}, 9(1990), 1--26. \re{[7]} Huishi Li, {\it
Noncommutative Gr\"obner Bases and Filtered-Graded Transfer}, LNM,
1795, Springer-Verlag, 2002. \re{[8]} Huishi Li, The general PBW
property, {\it Algebra Colloquium}, in press.
\re{[9]} T.~Mora, An introduction to commutative and noncommutative
Gr\"obner bases, {\it Theoretical Computer Science}, 134(1994),
131--173. \re{[10]} J.C.~McConnell and J.C.~Robson, {\it
Noncommutative Noetherian Rings}, John Wiley \& Sons, 1987.
\re{[11]} I. Peeva, V. Reiner and B. Sturmfels, How to shell a
monoid, {\it Math. Ann.}, 310(1998), 379¨C393.

\end{document}